\documentclass[12pt]{amsart}
\usepackage{amscd,amssymb}
\usepackage[graph,frame,poly,arc]{xy}  
\usepackage[plainpages,backref,urlcolor=blue]{hyperref}

\theoremstyle{plain}
\newtheorem{thm}[subsection]{Theorem}
\newtheorem{lem}[subsection]{Lemma}
\newtheorem{prop}[subsection]{Proposition}
\newtheorem{cor}[subsection]{Corollary}

\theoremstyle{definition}
\newtheorem{rk}[subsection]{Remark}
\newtheorem{definition}[subsection]{Definition}
\newtheorem{ex}[subsection]{Example}

\numberwithin{equation}{section}
\newcommand{\OO}{{\mathcal O}}

\newcommand{\K}{{\mathcal K}}

\newcommand{\SSS}{{\mathcal S}}

\newcommand{\C}{\mathbb{C}}
\newcommand{\PP}{\mathbb{P}}

\DeclareMathOperator{\codim}{codim}

\begin{document}

\title [Hessian ideals  and generalized Tjurina algebras]
{ Hessian ideals of a homogeneous polynomial and generalized Tjurina algebras}

\author[Alexandru Dimca]{Alexandru Dimca$^1$}
\address{Univ. Nice Sophia Antipolis, CNRS,  LJAD, UMR 7351, 06100 Nice, France. }
\email{dimca@unice.fr}

\author[Gabriel Sticlaru]{Gabriel Sticlaru}
\address{Faculty of Mathematics and Informatics,
Ovidius University,
Bd. Mamaia 124, 900527 Constanta,
Romania}
\email{gabrielsticlaru@yahoo.com }
\thanks{$^1$ Partially supported by Institut Universitaire de France and IAS Princeton.}



\begin{abstract} Using the minors in  Hessian matrices, we introduce new graded algebras associated to a homogeneous polynomial. When the associated projective hypersurface has isolated singularities, these algebras are related to some new local algebras associated to isolated hypersurface singularities, which generalize their Tjurina algebras.
One consequence of our results is a new way to determine the number of weighted homogeneous
singularities of such a hypersurface. 

\end{abstract}

\maketitle


\section{Introduction and statement of results} \label{sec:intro}

Let $S=\oplus_k S_k=\C[x_0,...,x_n]$ be the graded polynomial ring in $n+1$ indeterminates with complex coefficients, where $S_k$ denotes the vector space of degree $k$ homogeneous polynomials. Consider  for a polynomial $f \in S_d$, the corresponding Jacobian ideal $J_f$ generated by the partial derivatives $f_j$ of $f$ with respect to $x_j$ for $j=0,...,n$, the graded Milnor algebra $M(f)=\oplus_k M(f)_k=S/J_f$ and its Hilbert series
\begin{equation}
\label{HPMf}
 HP(M(f);t)= \sum_k\dim M(f)_kt^k.
\end{equation}
Assume in this note that the projective hypersurface $V=V(f):f=0$  in $\PP^n$ is  reduced  with (at most) isolated singularities at the points $a_1,...,a_p$. Let $Hess(f)$ be the Hessian matrix $(f_{ij})$ of the second order partial derivatives of $f$ and  $h(f)$ be the Hessian of $f$, i.e. the determinant of this matrix $Hess(f)$.
More generally, for each $k$ satisfying $1 \leq k \leq n+1$ we denote by $h_k(f)$ the ideal in $S$ generated by all $k \times k$-minors in the matrix $Hess(f)$. In particular, the ideal $h_{n+1}(f)=(h(f))$ is a principal ideal.
For each $k$ as above, consider the {\it graded $k$-th Hessian algebra of the polynomial $f$ } defined by
\begin{equation}
\label{def1}
H_k(f)=S/(J_f+h_k(f)),
\end{equation}
whose isomorphism class is clearly a $GL(n+1,\C)$-invariant of $f$. Clearly we have a sequence of epimorphisms of graded $\C$-algebras
$$M(f) \to H_{n+1}(f) \to H_n(f) \to \cdots \to H_1(f).$$
On the other hand, we introduce for any isolated hypersurface singularity $(V,0)$ at the origin of $\C^n$, given by an analytic germ $g=0$, the  $k$-th Hessian ideal $h_k(g)$ to be the ideal generated by  all the $k$-th minors in the Hessian $(n \times n)$-matrix $Hess(g)=(g_{ij})$. By convention, we set $h_{n+1}(g)=0.$ It is easy to check that the isomorphism class of the {\it local $k$-th Hessian algebra of the germ $g$ } defined by
\begin{equation}
\label{def2}
H_k(g)=\OO_n/((g)+J_g+h_k(g)),
\end{equation}
 is a $\K$-invariant of $g$, i.e. depends only on the isomorphism class of the germ $(V,0)$, see Lemma \ref{lem1}. Clearly we have now a sequence of epimorphisms of local Artinian $\C$-algebras
$$T(g) = H_{n+1}(g) \to H_n(g) \to \cdots \to H_1(g),$$
starting with the Tjurina algebra $T(g)$ of $g$. Hence one can consider the Hessian algebras $H_k(g)$ as a generalization of the Tjurina algebra $T(g)$.
For $k=1,...,n+1$  we introduce the $k$-th {\it Hessian number} of $g$ to be
\begin{equation}
\label{eq5} 
\chi_k(g)=\dim H_k(g).
\end{equation}
In view of the above remark, we can also write $\chi_k(g)=\chi_k(V,0)$ and use any normal form of $g$ to compute these new invariants. Note that $\chi_{n+1}(V,0)=\tau(V,0)$, the usual Tjurina number of $(V,0)$ and one has the inequalities
$$0\leq \chi_1(V,0) \leq ...\leq \chi_n(V,0) \leq \chi_{n+1}(V,0)=\tau(V,0).$$ 

Our first main result is the following.
\begin{thm}
\label{thm1}
Assume that $V(f)$ is a hypersurface  in $\PP^n$ with isolated singularities at the points $a_1,...,a_p$. Then for $k=1,2,...,n+1$ and any $m$ large enough one has
$$\dim H_k(f)_m=\sum_{i=1,p} \chi_k(V(f),{a_i}).$$

\end{thm}

The following result estimates how large $m$ should be in order to have such stabilization results.
First we recall some definitions, see \cite{DStEdin}.

\begin{definition}
\label{def}

For a hypersurface $V:f=0$ with isolated singularities we introduce three integers, as follows.

\noindent (i) the {\it coincidence threshold} 
$$ct(V)=\max \{q:\dim M(f)_k=\dim M(f_s)_k \text{ for all } k \leq q\},$$
with $f_s$  a homogeneous polynomial in $S$ of degree $d$ such that $V_s:f_s=0$ is a smooth hypersurface in $\PP^n$.

\noindent (ii) the {\it stability threshold} 
$st(V)=\min \{q~~:~~\dim M(f)_k=\tau(V) \text{ for all } k \geq q\},$
where $\tau(V)$ is the total Tjurina number of $V$, that is
$\tau(V)=\sum_{i=1,p} \tau(V(f),{a_i}).$

\noindent (iii) the {\it minimal degree of a  syzygy} $mdr(V)=\min \{q~~:~~ H^n(K^*(f))_{q+n}\ne 0\}$,
where $K^*(f)$ is the Koszul complex of $f_0,...,f_n$ with the natural grading.

\end{definition}

It is known that one has
\begin{equation} 
\label{REL}
ct(V)=mdr(V)+d-2,
\end{equation} 
and precise estimates on $mdr(V)$ are given in \cite{DTor}, \cite{DS2}, \cite{DStNach}, \cite{DS14}. Our second main result is the following.

\begin{thm}
\label{thm2}
Assume that the hypersurface $V(f)$ in $\PP^n$ has  only isolated  singularities.
Then one has the following stabilization property: 
$$ \dim H_n(f)_m=\tau(V(f))-|Sing_{wh}(V(f))|$$
for any $m$ large, where $|Sing_{wh}(V(f))|$ is the number of weighted homogeneous singularities the hypersurface $V(f)$ has.
Moreover, for any integer $k$, $1 \leq k\leq n$, the dimension of $H_k(f)_m$ is constant for any 
$$m \geq T_k:=\max(T-ct(V)+k(d-2), st(V)). $$
\end{thm}

We also have the following result, expressing $HP(H_{n+1}(f);t)$ in term of the series $HP(M(f);t)$, which was already studied in \cite{DS2}, \cite{DStNach}, \cite{DStEdin}.

\begin{prop}
\label{propA}
(i)   Assume the hypersurface $V(f)$ is smooth. Then
$$HP(H_{n+1}(f);t)=HP(M(f);t)-1=\frac{(1-t^{d-1})^{n+1}}{(1-t)^{n+1}}-1.$$

(ii) Assume the hypersurface $V(f)$ has isolated singularities (and it is not smooth). Then
$$HP(H_{n+1}(f);t)=HP(M(f);t).$$
\end{prop}

A consequence of these results is the following (possibly the fastest) way to compute the total number of singularities of a projective hypersurface in cases when it is known a priori that these singularities are all weighted homogeneous.

\begin{cor}
\label{corE}
Assume that the hypersurface $V(f)$ in $\PP^n$ has  only isolated  singularities (and it is not smooth).
Then the total number of weighted homogeneous singularities of $V(f)$ is given by
$$ |Sing_{wh}(V(f))|=\dim M(f)_m  - \dim H_n(f)_m$$
for any $m\geq \tilde T_n=(2n+1)(d-2)$.
\end{cor}

In the second section we prove these results, give some Corollaries and a number of Examples.
In the final section we offer a discussion on the constructible partitions induced on the space of homogeneous polynomials of a fixed degree $d$ by the values of the Hilbert-Poincar\'e series of the Milnor algebra $M(f)$ and of the Hessian algebras $H_k(f)$. The only results here are of an {\it experimental computational nature}. The computations of various invariants given in this paper were made using two computer algebra systems, namely CoCoA \cite{Co} and Singular \cite{Sing}.
The corresponding codes are available on request.

\section{Hessian matrix and Hessian ideals} 

\begin{lem}
\label{lem1}
The isomorphism class of the  local $k$-th Hessian algebra 
$$H_k(g)=\OO_n/((g)+J_g+h_k(g)),$$
of the isolated hypersurface germ $(V,0):g=0$  is a $\K$-invariant of $g$ for any integer $k=1,2,...,n$.

\end{lem}

\proof
We have to show that if $u \in \OO_n$ is a unit and $\phi: (\C^n,0) \to (\C^n,0)$ an isomorphism germ, and if we set $g'(y)=u(y) \cdot g(\phi(y))$, then the algebra $H_k(g')$ is isomorphic to the algebra $H_k(g)$ for any $k$. Note first that the quotients $\OO_V=\OO_n/(g)$ and $\OO_{V'}=\OO_n/(g')$ are isomorphic $\C$-algebras, since the isomorphism $\phi^*: \OO_n \to \OO_n$ given by $a(y) \mapsto a(\phi(y))$ sends the ideal $(g)$ injectively onto the ideal $(g')$. Computing the partial derivatives of $g'$ with respect to $y_1,...,y_n$, we see that the isomorphism $\phi^*$ sends also the ideal $(g)+J_g$ into the ideal $(g')+J_{g'}$.
The inclusion $\phi^*((g)+J_{g}) \subset (g')+J_{g'}$ implies $\tau(g) \geq \tau(g')$. Since the relation between $g$ and $g'$ is symmetrical, we get in the same way  $\tau(g') \geq \tau(g)$, and hence in fact we have equalities everywhere.
Now compute the second order partial derivatives of $g'$ and note that, modulo the ideal $(g')+J_{g'}$, the Hessian matrix $Hess(g')$  is a product of 3 matrices $A(y)^t Hess(g)(\phi (y))A(y)$, where $A(y)=d\phi(y)$ is the Jacobian matrix of $\phi$ and $A(y)^t$ denotes its transpose. From this relation it follows that 
$$\phi^* ((g) +J_g+h_k(g)) \subset (g')+J_{g'}+m_k(g')$$
for any $k$. The equality is established exactly as above, using the symmetry between $g$ and $g'$.

\endproof

\subsection{Proof of Theorem \ref{thm1}.} 

For any $k$, it is clear that the ideal $J_f+h_k(f)$ defines a 0-dimensional subscheme in $\PP^n$ whose support is contained in the singular locus $V(f)_{sing}=\{a_1,...,a_p\}$. Hence one has for large $m$ the following equality
\begin{equation}
\label{eq1}
\dim H_k(f)_m=\sum_{i=1,p} \dim H_k(f)_{a_i},
\end{equation}
where $H_k(f)_{a_i}$ denotes the analytic localization of the graded $S$-module $H_k(f)$ at the point $a_i$, exactly as in the proof of Corollary 9 in \cite{CD}.

To study such a localization, assume that the coordinates have been chosen such that the singular point $a=a_1$ is located at $(1:0:...:0)$. Then the isolated hypersurface singularity $(V(f),a)$ has a local equation 
$$g(y_1,...,y_n)=f(1,y_1,...,y_n)=0.$$
It follows that $g_j=f_j(1,y_1,...,y_n)$ for any $j=1,...,n$, where $g_j$ denotes the partial derivative of $g$ with respect to $y_j$. Using the Euler relation for $f$ exactly as in the proof of Corollary 9 in \cite{CD}, we also get
\begin{equation}
\label{eq2}
d\cdot g=f_0(1,y_1,...,y_n)+y_1f_1(1,y_1,...,y_n)+...+y_nf_n(1,y_1,...,y_n).
\end{equation}
This implies that the localization of the Jacobian ideal $J_f$ at the point $a$ coincide with the ideal spanned by $g, g_1,...,g_n$, i.e. the Tjurina ideal of the isolated hypersurface singularity $(V(f),a)$.
Let's look now at the second order derivatives. One clearly has $g_{ij}=f_{ij}(1,y_1,...,y_n)$ for all
$1 \leq i,j \leq n$. If we take the partial derivative with respect to $y_j$ for some $j>0$ of the polynomials entering into the equality \eqref{eq2}, we get
\begin{equation}
\label{eq3}
d\cdot g_j=f_{0j}(1,y_1,...,y_n)+y_1f_{1j}(1,y_1,...,y_n)+...+
\end{equation}
$$+(f_j(1,y_1,...,y_n)+y_jf_{jj}(1,y_1,...,y_n))+...+y_nf_{nj}(1,y_1,...,y_n).$$
Similarly, starting with the Euler relation for $f_0$, taking the partial derivative with respect to $x_0$ and then setting $x_0=1$ and $x_j=y_j$ for $j>0$ we get
\begin{equation}
\label{eq4}
(d-1)\cdot f_0(1,y_1,...,y_n)=(f_0(1,y_1,...,y_n)+f_{00}(1,y_1,...,y_n))+ \sum_{j>0}y_j
f_{0j}(1,y_1,...,y_n).
\end{equation}
The formulas \eqref{eq3} and \eqref{eq4} give a proof of the following result.

\begin{prop}
\label{prop1}
The localization of the ideal $J_f+h_k(f)$ at the singular point point $a=(1:0:...:0)$ coincides with the ideal in the corresponding local ring $\OO_n$ spanned by the local equation $g(y)=f(1,y_1,...,y_n)$, its partial derivatives $g_j$ with respect to $y_j$ for $j=1,...,n$ and all $k\times k$ minors in the local $n \times n$ Hessian matrix $Hess(g)=(g_{ij})$.
\end{prop}

\proof Indeed, modulo the Tjurina ideal, the matrix $Hess(f)(1,y_1,...,_n)$ can be transformed using the above formulas by obvious line and column operations into a matrix obtained from the local Hessian matrix $Hess(g)$ by adding a first line and a first column formed by zeros only.

\endproof

Combining this definition with Proposition \ref{prop1} and the equality \eqref{eq1} clearly completes the proof 
of Theorem \ref{thm1}.

\begin{ex}
\label{exnodal}
Assume that $(V,0)$ is an $A_1$-singularity, hence a node. Then we can choose an equation of the form $g=y_1^2+...+y_n^2$ and it is clear that $\chi_k(V,0)=0$ for $k=1,2,...,,n$ and $\chi_{n+1}(V,0)=1$. In particular, if $V(f)$ is a nodal hypersurface then $ H_k(f)_m=0$ for $1\leq k \leq n$ and hence the corresponding Hilbert-Poincar\'e series
$$HP(H_k(f);t) =\sum_m \dim H_k(f)_m t^m$$
is a polynomial.
\end{ex}

\begin{ex}
\label{excusp}
Assume that $n=2$ and $(V,0)$ is an $A_q$-singularity. Then we can choose an equation of the form $g=y_1^2+y_2^{q+1}$ and it is clear that 
$\chi_1(V,0)=0$,  $\chi_2(V,0)=q-1$ and $\chi_3(V,0)=q$. In particular, if $V(f)$ is a cuspidal curve then 
$  \dim H_k(f)_m=(k-1)\kappa$, where $\kappa$ is the number of cusps $A_2$.
To have a numerical example, let $f=x^2y^2+y^2z^2+x^2z^2-2xyz(x+y+z)$. Then $f=0$ has $\kappa=3$ cusps and no other singularities, see also \cite{D1}, p.129 for more on this quartic curve.
Numerical computations  show that
$$HP(H_1(f);t)=1+3t=1+3t+0(t^2+t^3+...)$$
$$HP(H_2(f);t)=1+3t+6t^2+7t^3+3(t^4+t^5+...)$$
$$HP(H_3(f);t)=1+3t+6t^2+7t^3+6(t^4+t^5+...).$$
In particular
$$ \dim H_k(f)_m=3(k-1)$$
for $k=1,2,3$ and $m\geq 4$.

\end{ex}

\begin{ex}
\label{exD}
Assume that $n=2$ and $(V,0)$ is an $D_q$-singularity, $q \geq 4$. Then we can choose an equation of the form $g=y_1^2y_2+y_2^{q-1}$ and it is clear that 
$\chi_1(V,0)=1$,  $\chi_2(V,0)=q-1$ and $\chi_3(V,0)=q$.
Moreover, it clear that for a plane curve singularity $(V,0)$, one has $\chi_1(V,0)=1$, if and only if $(V,0)$ has type $D_q$ for some $q \geq 4$.

\end{ex}

Using Theorem \ref{thm1} and Example \ref{excusp} we get the following.

\begin{ex}
\label{exC}
The curve $V(f)$ in $\PP^2$ has  only $A_q$ singularities  for various $q$'s (i.e. $V(f)$ has only corank 1 singularities) if and only if  $  H_1(f)_m=0$ for $m$ large.
If this holds, then the total number of singularities of $V(f)$ is given by
$$ |Sing(V(f))|=\dim H_3(f)_m  - \dim H_2(f)_m $$
for $m$ large.
\end{ex}

\subsection{Proof of Theorem \ref{thm2}.}

Consider the local case  of an isolated weighted homogeneous singularity $(V,0)$ given by $g=0$.
Then $g \in J_g$ and in the local Gorenstein ring $M(g)=\OO_n/J_g$, the Hessian $h(g)$ of $g$ spans the socle which is one dimensional, see \cite{EGH}. Clearly we have $H_n(g)=M(g)/(h(g))$, and hence $\chi_n(g)=\tau(g)-1=\chi_{n+1}(g)-1$. Next, if $g$ is not weighted homogeneous, then the ideal $(g)$ in $M(g)$ is nonzero, and hence contains the minimal ideal $(h(g))$. It follows that
in such a case $\chi_n(g)=\tau(g)$.
The result then follows  for large $m$ using 
Theorem \ref{thm1} . 

Let $I$ be the saturation of the Jacobian  ideal $J_f$. Then Prop. 2. in \cite{DBull} tells us that
\begin{equation}
\label{eq10} 
\dim S_q/I_q=\tau(V(f))
\end{equation}
for any $q \geq T-ct(D)$. In other words, for such $q$'s, the evaluation map
\begin{equation}
\label{eq10.5} 
ev_q: S_q/I_q \to \oplus \OO_n/((g)+J_g)
\end{equation}
is a bijection, where the sum is taken over all singularities $(V,a):g=0$ of the hypersurface $V=V(f)$. In particular, for any fixed singularity $a \in V_{sing}$ there is a polynomial $h_a \in S_q$ such that the class of $h_a$ in the summand $T(g)=\OO_n/((g)+J_g)$ corresponding to $a$ is $1$ and its classes in all the other summands are $0$. The $k$-minors in the matrix $Hess(f)$ have degree $k(d-2)$ and the computation done in Proposition \ref{prop1} shows that for any singular point $a \in V_{sing}$ there are linear combination $m_s$ of such minors (arising from a change of coordinates on $\PP^n$) such that the classes of $m_s$ in the corresponding local ring $T(g)$
generate the ideal $h_p(g)$.
Cor. 2 in \cite{DBull} tells us that
\begin{equation}
\label{eq11} 
I_q=J_{f,q}
\end{equation}
for all $q \geq q_1=\max (T-ct(V(f)), st(V(f)))$.

Consider the bijection $ev_q$ for $q \geq T_k:=\max(T-ct(V)+k(d-2), st(V)) \geq q_1$. Since all the products $h_am_s \in h_k(f)_q$, it follows that 
$$ev_q((J_{f,q}+h_k(f)_q)/J_{f,q})=\oplus ((g)+J_g+h_k(g)/((g)+J_g).$$
This gives the following bijection for $q \geq T_k$
$$ev_q: H_k(f)_q \to \oplus \OO_n/((g)+J_g+h_k(g))$$
which completes the proof of Theorem \ref{thm2}.

\begin{rk}
\label{rkbound}
Note that $T_k \leq \tilde T_k=T+k(d-2)=(n+k+1)(d-2)$ for any $k$. So one can use this larger bound when there is no information on $ct(V)$ and $st(V)$, i.e. we always have stabilization for $\dim H_k(f)_m$ when $m \geq \tilde T_k$. When we have information on $ct(V)$ but not on $st(V)$, as it is often the case, one can use the better general bound
$$\widehat T_k=T+\max(k(d-1)-ct(V),0).$$

\end{rk}

\begin{cor}
\label{corA}
(i) The hypersurface $V(f)$ in $\PP^n$ is  nodal if and only if  $ H_n(f)_m=0$ for 
$m\geq \widehat  T_n=3n(d-2)/2.$ In particular, for $n=2$, we have $H_2(f)_m=0$ for $m\geq T= \widehat T_2=3d-6$.

(ii) Assume the hypersurface $V(f)$ in $\PP^n$ is  nodal and $n \geq 4$. Then $H_{n-1}(f)_m=0$ for 
$m\geq \widehat  T_{n-1}=(3n-2)(d-2)/2.$
\end{cor}
\proof

For a nodal hypersurface, it is known that $ct(V(f)) \geq (n+2)(d-2)/2$, see \cite{DS2}.
It follows that 
$$T-ct(V(f))+n(d-2) \leq T+(n-2)(d-2)/2=3n(d-2)/2 \geq T,$$ 
 which proves the result (i). The second claim (ii) follows along the same lines. The condition $n\geq 4$ is necessary to insure that $\widehat T_{n-1} \geq T$.

\endproof

\begin{ex}
\label{nodalcurves}
For the nodal curves $f=xyz^{d-2}+x^d+y^d=0$, a direct computation  shows that the stabilization $H_2(f)_m=0$ occurs for $m \geq 3d-7$ for low values of $d$, so only one unit better than our prediction. For the Chebyshev curves considered in \cite{DStCamb} or in \cite{St}, a direct computation  shows that the stabilization $H_2(f)_m=0$ occurs for $m \geq 3d-8$ for low values of $d$, so just 2 units better than our general prediction in Corollary \ref{corA}.
\end{ex}

\begin{ex}
\label{nodal4folds}
Let $V:f=0$ be the Chebyshev 4-fold of degree $d$ in $\PP^5$ considered in \cite{DStNach} or in \cite{St}. A direct computation  shows that the stabilization $H_4(f)_k=0$ occurs for 
$k \geq 13m-10$ for low values of $d=2m+1$, and for $k \geq 13m-17$ for low values of $d=2m$.
The corresponding (integral) values given by $T'_4$ in these cases are  $13m-7$ for low values of $d=2m+1$ and $13m-13$ for low values of $d=2m$, so the difference between the theoretical bound and the real value of $k$ is less than 4 for any degree $d$.
Similarly, the stabilization $H_5(f)_k=0$ occurs for 
$k \geq 15m-10$ for low values of $d=2m+1$, and for $k \geq 15m-17$ for low values of $d=2m$.
The corresponding (integral) values given by $T'_5$ in these cases are  $15m-8$ for low values of $d=2m+1$ and $15m-15$ for low values of $d=2m$, so the difference between the theoretical bound and the real value of $k$ is less than 2 for any degree $d$.
\end{ex}

\begin{cor}
\label{corB}
Assume the curve $V(f)$ in $\PP^2$ has $\nu$ nodes and $\kappa$ cusps as singularities.
Then $  \dim H_2(f)_m =\kappa$ and $  \dim M(f)_m=\dim H_3(f)_m =\nu +2\kappa$ for $m\geq T+d/6$. In other words, the asymptotic behavior of $\dim H_k(f)_m$ for $k=2,3$ determines the number of nodes and cusps.
\end{cor}

\proof

For a plane curve with nodes and cusps, it is known that $ct(V(f)) \geq 11d/6-4$, see \cite{DS14}.
It follows that 
 $T_1 \leq T+ d/6$, which proves the result.

\endproof

The stabilization value of $H_{n-1}(f)_m$ is not easy to describe gemetrically beyond Theorem \ref{thm1}, as the following example shows.

\begin{ex}
\label{simplecurves}
Let $V:f=0$ be a plane curve having only simple singularities, hence singularities of type $A_k$, $D_k$, $E_6$, $E_7$ and $E_8$. A direct computation shows that $\chi_1(g)=2$ for $g$ of type $E_6$ and $E_7$, and  $\chi_1(g)=3$ for $g$ of type $E_8$. Using now Example \ref{exD}, it follows that, for $m \geq \tilde T_1=4d-8$, one has 
$$\dim H_1(f)_m= \sum_k\sharp D_k+2(\sharp E_6 +\sharp E_7)+3 \sharp E_8.$$
\end{ex}
 
\begin{ex}
\label{exlines}
Consider the line arrangement in $\PP^2$ given by the following equation
$$f=(x^2-y^2)(y^2-z^2)(x^2-z^2).$$
Then the curve $V(f):f=0$ is a union of 6 lines and it has 3 nodes $A_1$ and 4 triple points $D_4$.
In this case $ct(V)=6$ as shown by $HP(H_3(f);t)=HP(M(f);t)$,  and hence $T_2=\max(12-6+8,st(V))=14$ and $T_1=T_2-(d-2)=10.$
Numerical computations  give the following.
$$HP(H_1(f);t)=1+3t+6t^2+10t^3+10t^4+7t^5+5t^6+4(t^7+t^8+...)$$
(hence again a stabilization 3 units before $T_1$)
$$HP(H_2(f);t)= 1+3t+6t^2+10t^3+15t^4+18t^5+19t^6+19t^7+13t^8+12(t^9+t^{10}+...)$$
(hence again a stabilization 2 units before $T_2$) and 
$$HP(H_3(f);t)=HP(M(f);t)=1+3t+6t^2+10t^3+15t^4+18t^5+19(t^6+t^7+...),$$
as predicted by Proposition \ref{propA}.

\end{ex}

The following consequence of the proof of Theorem \ref{thm2} might be useful sometimes.

\begin{cor}
\label{corD}
An isolated hypersurface singularity $g=0$ is weighted homogeneous if and only if $\chi_n(g)=\tau(g)-1$. Otherwise $\chi_n(g)=\tau(g)$. 
\end{cor}

\subsection{Proof of Proposition \ref{propA}}

In the case (i), it is enough to note that $M(f)$ is Artinian Gorenstein ring whose 1-dimensional socle is spanned by Hessian $h(f)$ which has degree $T$, see for instance \cite{EGH}.
\endproof

In the case (ii), the class of the Hessian $h(f)$ will now vanish in the local Tjurina ring $T(g)$ attached to any singularity of $V(f)$ by Proposition \ref{prop1}.  Indeed, recall that the first row and first column are zero in the Hessian matrix $Hess(f)$ (in suitable coordinates) modulo the Tjurina ideal $(g)+J_g$. Moreover, as in the proof of Theorem \ref{thm2}, the evaluation map
\begin{equation}
\label{eq20} 
ev_T: S_T/J_{f,T} \to \oplus \OO_n/((g)+J_g)
\end{equation}
is an isomorphism. It follows that $h(f) \in J_f$ and hence $H_{n+1}(f);t)=HP(M(f);t)$ in this case.
\endproof

The Hilbert-Poincar\'e series of the Hessian rings $H_k(f)$ for $1\leq k \leq n$ are quite mysterious, as the following example shows.

\begin{ex}
\label{2smooth}
Consider the Fermat curve $V(f): f=x^6+y^6+z^6=0$ and its "deformation"
$V(f'): f'=f+2xyz^4+3xzy^4+5yzx^4=0.$ One obtains following results by direct computation.
$$HP(M(f);t)=HP(M(f');t)=1+3t+6t^2+10t^3+15t^4+18t^5+19t^6+18t^7+15t^8+$$
$$+10t^9+6t^{10}+3t^{11}+t^{12}.$$
In particular, $V(f')$ is smooth as well and   $HP(H_3(f);t)=HP(H_3(f');t)$ can be computed using   Proposition \ref{propA}. For $k=2$ we get the following results.

\medskip

\noindent $HP(H_2(f);t)=1+3t+6t^2+10t^3+15t^4+18t^5+19t^6+18t^7+12t^8+7t^9+3t^{10}$

\medskip

\noindent and

\medskip

 \noindent  $HP(H_2(f');t)=1+3t+6t^2+10t^3+15t^4+18t^5+19t^6+18t^7+9t^8.$

\medskip

\noindent Note that the 2-minors have degree 8, which explains the fact that the monomials containing $t^m$ for $m<8$ are unchanged. For $f$ there are only 3 such minors, which explain the drop from 15 to 12 in the coefficient of $t^8$. For $f'$, there are 6 distinct minors (due to the symmetry of the Hessian matrix) and these 2-minors span a 6-dimensional vector space modulo $J_{f',8}$ in this case. Finally, for $k=1$ we get

\medskip

\noindent $HP(H_1(f);t)=1+3t+6t^2+10t^3+12t^4+12t^5+10t^6+6t^7+3t^8+t^9$

\medskip

\noindent and

\medskip

 \noindent  $HP(H_1(f');t)=1+3t+6t^2+10t^3+9t^4+3t^5.$

\medskip

\noindent Note that the 1-minors have degree 4, which explains the fact that the monomials containing $t^m$ for $m<3$ are unchanged. For $f$ there are only 3 such minors, which explain the drop from 15 to 12 in the coefficient of $t^4$. For $f'$, there are 6 distinct 1-minors (due to the symmetry of the Hessian matrix) and these 1-minors are linearly independent  in this case
as seen by looking at the coefficient of $t^4$. Moreover, if we look at the coefficient of $t^5$ in 
$HP(H_1(f');t)$, it says that {\it the only relations among the 18 polynomials $x_mf'_{ij}$ with $i\leq j$ and the 3 partial derivatives $f'_0,f'_1,f'_3$ are the obviuos 3 relations given by the Euler's formulas}
$$(d-1)f'_q=x_0f'_{q0}+x_1f'_{q1}+x_2f'_{q2},$$
for $q=0,1,2$. This accounts for the drop $18-15=3$ when we pass from $M(f')$ to $H_1(f')$.

\noindent The symmetry of the polynomial $HP(H_1(f);t)$ is easy to explain in general. Indeed, if we denote by $F(d)$ the Fermat polynomial in $n+1$ variables of degree $d$, then it is clear that 
$$H_1(F(d))=M(F(d-1)).$$

\end{ex}

\section{Constructible partitions for the space of homogeneous polynomials} \label{sec2}

 Consider the projective space $\PP(S_d)$
of non-zero homogeneous polynomials of degree $d$ and let $G=PGL_n(\C)$ act on $\PP(S_d)$ in the usual way, i.e. by linear changes of coordinates. Then there is a Zariski open and dense, $G$-invariant subset $U_s(n,d)$ in $\PP(S_d)$, parametrizing the smooth hypersurfaces of degree $d$ in $\PP^n$. There is a larger $G$-invariant Zariski open subset $U_{iso}(n,d)$ such that $f \in U_{iso}(n,d)$ if and only if the hypersurface $V(f):f=0$ in $\PP^n$ has only isolated singularities.
Moreover, for each list of (distinct or not) simple isolated hypersurface singularities $(Y_1,Y_2,...,Y_p)$ we can consider the constructible
subset
$S(Y_1,Y_2,...,Y_p) \subset  U_{iso}(n,d)$ defined by the condition that $f \in S(Y_1,Y_2,...,Y_p)$ if and only if the hypersurface $V(f)$ has exactly $p$ singularities of types $Y_1,...,Y_p$.
It is known that if this stratum is non-empty, then
$$\codim S(Y_1,Y_2,...,Y_p) \leq \sum _{i=1,p}\tau(Y_i)$$
with equality if the degree $d$ is large, e.g. if
\begin{equation}
\label{connect}
d \geq \sum _{i=1,p}s_i+p-1, 
\end{equation}
where $s_i$ is the $\K$-determinacy order of $Y_i$, see Prop. 1.3.9 in \cite{D1}, p.18. Moreover in this case the stratum $S(Y_1,Y_2,...,Y_p)$ is smooth and connected. For more precise and more general results, see Theorem 1 and Theorem 5 in \cite{GLS}.

When all of the singularities $Y_i$ are nodes and $n=2$, then the corresponding stratum $S_0(Y_1,Y_2,...,Y_p)$ consisting of all reduced {\it irreducible}  curves with $p$ nodes is connected. This is the famous Severi problem solved by J. Harris \cite{Harrris}. However, in general, the strata 
 $S(Y_1,Y_2,...,Y_p)$ are not connected even in this case, see Example \ref{ex3A1} below. For more general and precise results in this direction, as well as for an extensive bibliography on the subject, we refer to the excellent survey \cite{GLS}.

Since $\dim M(f)_k= \tau (V(f))$, for $k >T=(n+1)(d-2)$, see \cite{CD}, it follows that there is a finite constructible {\it Jacobian partition} $\SSS_Y(M(f)$ of each stratum $S(Y)=S(Y_1,Y_2,...,Y_p)$, such that $f, g \in S(Y)$ are in the same stratum of $\SSS_Y$ if and only if $HP(M(f);t)=HP(M(g);t)$.
It is known that $f_s \in U_s(n,d)$ if and only if
\begin{equation}
\label{eqFermat} 
HP(M(f_s);t)=\frac{(1-t^{d-1})^{n+1}}{(1-t)^{n+1}},
\end{equation}
and hence the Jacobian partion on $U_s(n,d)$ consists of just one stratum.
However, for a nodal hypersurface $V(f)$ with 3 nodes, the value of $HP(M(f);t)$ depends on whether the nodes are collinear or not, see Example 4.3 (iii) in \cite{DStEdin}. It follows that the Jacobian stratification on the stratum $S(3A_1)$ is not trivial in general.

\begin{ex}
\label{ex3A1}
Consider the case $n=2$,  $d=4$ and the stratum $S(3A_1)$. It is known that this stratum 
$S(A_3)$ has two irreducible components, $X=S_0(3A_1)$ containing the irreducible quartics with 3 nodes (necessarily non-collinear) and $Y$, the reduced quartics with 3 collinear nodes obtained from a smooth cubic and a trisecant line, see Remark 4.7.21 in \cite{Se}. Then one has $\codim X= \codim Y = 3$, as follows from Theorem 4.7.18 in \cite{Se}.

  Example 4.3 (iii) in \cite{DStEdin} tells us that there are two possibilities for $H(M(f);t)$ on the stratum $S(3A_1)$, namely
the {\it minimal value} (coefficient by coefficient) is 
$$HP(M(f);t) = 1+3t + 6t^2 + 7t^3 + 6t^4 + 3t^5 + 3t^6+ ...$$
and this is attained exactly on $X$ and the other possible value is
$$HP(M(g);t) = 1+3t + 6t^2 + 7t^3 + 6t^4 + 4t^5 + 3t^6+ ...$$
and which is attained precisely on $Y$. In other words, in this case the Jacobian partition separates the irreducible components of  the stratum $S(3A_1)$.

Consider now the case $n=2$, but $d \geq 8$. Since the $\K$-determinacy order of a node is $s_i=2$, it follows from \eqref{connect} that the stratum $S(3A_1)$ is smooth and connected in this case. The minimal value of $HP(M(f);t)$, which is
\begin{equation}
\label{HPmin} 
\frac{(1-t^{d-1})^{3}}{(1-t)^{3}}+ 2t^T+3t^{T+1}(1+t+t^2+...),
\end{equation}
by Example 4.3 (iii) in \cite{DStEdin}, is attained on a Zariski open dense subset $S'$ of the stratum $S(3A_1)$, while the larger value
\begin{equation}
\label{HPlarge} 
\frac{(1-t^{d-1})^{3}}{(1-t)^{3}}+ t^{T-1}+2t^T+3t^{T+1}(1+t+t^2+...),
\end{equation}
is attained on the closed subset $S''=S(3A_1) \setminus S'$ in $S(3A_1)$.
\end{ex}

\begin{ex}
\label{ex6A2}
Consider the case $n=2$,  $d=6$ and the stratum $S(6A_2)$. It is known by the classical work of Zariski that this stratum 
 has two irreducible components, $X$ containing the irreducible sextics with 6 cusps on a conic and $Y$, containing the irreducible sextics with 6 cusps not on a conic.
As an example of a sextic in $X$ one can take 
$$V(f): f=(x^2+y^2)^3+(y^3+z^3)^2=0$$
and then a direct computation  yields 
$$ HP(M(f);t) = 1+3t+6t^{2}+10t^{3}+15t^{4}+18t^{5}+19t^{6}+18t^{7}+16t^{8}+13t^{9}+12(t^{10}+.... .$$
As an example of a sextic in $Y$ one can take 
$$V(g): g=27(x+y)^3(x+y-z)^2(x+y+2z)-27(x+y)^2(x+y-z)^2((x-y)^2-z^2)+$$
$$+9((x+y)^2-z^2)((x-y)^2-z^2)^2-((x-y)^2-z^2)^3=0,$$
see \cite{Oka}, formula (5.7), and again a direct computation  yields 
$$HP(M(g);t) = 1+3t+6t^{2}+10t^{3}+15t^{4}+18t^{5}+19t^{6}+18t^{7}+15t^{8}+12(t^{9}+....$$
We do not know whether the Hilbert-Poincar\'e series is constant on $X$ and/or $Y$, but it seems that again it can separate distinct irreducible components in a stratum.

\end{ex}

We have seen above, especially in Example \ref{2smooth}, that though one can also define {\it the Hessian partition} by looking at the strata inside a given $S(Y_1,...Y_p)$ where the corresponding Hilbert-Poincar\'e series $HP(H_k(f);t)$ are constant for $k=1,2,...,n+1$, the behavior of these partitions is rather mysterious. Indeed, Example \ref{2smooth} shows that already the Hessian partition on the stratum $U_s(n,d)$ can be quite subtle.

To investigate the possibilities for the Hilbert series for $HP(M(f);t)$ as well as for  $HP(H_k(f);t)$ for various $k$'s, in the neighborhood of a given point $g \in \PP(S_d)$, one can proceed as follows. It is known that the tangent space to the $G$-orbit of $g$ (where all the invariants are clearly preserved) is given by $(J_g)_d$. Let $h_1,...,h_q$ be a (monomial) basis of $M(g)_d$ and consider the family of polynomials
\begin{equation}
\label{unfold1}
f=g+a_1h_1+...a_qh_q,
\end{equation}
with $a_j \in \C$, which is exactly a transversal to the $G$-orbit of $g$. If we choose the $a_j$ small enough, the values obtained for $HP(M(f);t)$ as well as for  $HP(H_k(f);t)$ for various $k$'s, will cover all the possibilities that occur in a neighborhood of $g$ (in the strong topology). If we allow any values for $a_j$, they will give sometimes interesting values, even though not related to the fixed polynomial $g$ (but covering this time a Zariski open neighborhood of $g$).

\begin{rk}
\label{rkdeformation}
Consider the case of smooth plane curves of increasing degree $d$.

\noindent (i) When $d=2$, then $U_s(2,2)$ is a single $G$-orbit, hence all the invariants $HP(H_k(f))$ are constant on $U_s(2,2)$, i.e. the Hessian partition is trivial. This fact clearly holds  for any $n \geq 2$ as well.

\noindent (ii) When $d=3$, it is known that any polynomial in $U_s(2,3)$ is $G$-equivalent to a Hesse normal form
$$f=x^3+y^3+z^3+3axyz,$$
with $a^3 \ne -1$. This comes essentially from the fact that if $g=x^3+y^3+z^3$ is the corresponding Fermat polynomial, then $xyz$ is a basis for the 1-dimensional vector space $M(g)_3$.

\noindent (iii) When $d=4$, if we take $g=x^4+y^4+z^4$, then a basis for $M(g)_4$ is given by the following 6 monomials
$$x^2y^2, x^2z^2, y^2z^2,x^2yz,xy^2z,xyz^2.$$
But doing computation with 6 parameter families seems too complicated for the moment.
That is why in the following examples we take $h_1,...,h_q$  just {\it a family of independent elements} of $M(g)_d$, i.e. we explore what happens only in {\it a fixed direction in the normal space} to the orbit of $g$.

\end{rk}

\begin{ex}
\label{3para}
Consider the case $n=2$, $d=6$, $g=x^6+y^6+z^6$ and the 3-parameter family  
$$E:f=g+ ax^4yz+bxy^4z+cxyz^4,$$
which is a partial transversal as explained above.
We will describe all the possibilities for the series $HP(H_1(f);t)$ when $f$ has this form and $f \in U_s(2,6)$. By direct computation one gets the following Hessian strata in $E$, where for each strata we list the corresponding series $HP(H_1(f);t)$ and  the conditions on $a,b,c$ to get such a series, up to the action of the permutation group on $x,y,z$.

\medskip

\noindent $S_1: 1+3t+6t^2+10t^3+9t^4+3t^5  \text{ with } a,b,c \text{ generic }.$

\medskip

\noindent This is the largest stratum, a Zariski open dense subset in $E=\C^3$.

\medskip

\noindent $S_2: 1+3t+6t^2+10t^3+9t^4+4t^5  \text{ with } a=b \ne 0  \text{ and } c \text{ generic }.$

\medskip

\noindent $S_3: 1+3t+6t^2+10t^3+9t^4+6t^5  \text{ with } a=b = c \text{ generic }.$

\medskip

\noindent $S_4: 1+3t+6t^2+10t^3+9t^4+4t^5+t^6  \text{ with } a=b \ne 0   \text{ and }c =0.$

\medskip

\noindent $S_5: 1+3t+6t^2+10t^3+9t^4+6t^5+2t^6 \text{ with } a,b,c \text{ a collection of finite nonzero values},$ for instance $a=b=c=2$.

\medskip

\noindent $S_6: 1+3t+6t^2+10t^3+9t^4+7t^5 +3t^6+t^7 \text{ with } a=b =0 \text{ and } c \text{ generic }.$

\medskip

\noindent $S_7: 1+3t+6t^2+10t^3+12t^4+12t^5 +10t^6+6t^7+3t^8+t^9 \text{ with } a=b = c =0.$

\medskip

In the next picture (a Hasse diagram as in the theory of hyperplane arrangements) we draw a line between $S_i$ and $S_j$ with $S_i$ above $S_j$ if the 
inequality $HP(H_1(f^i);t)\geq HP(H_1(f^j);t)$ (coefficient by coefficient) holds for some $f^i \in S_i$ and $f^j \in S_j$ and we call the resulting poset a Hessian poset. This inequality is {\it a necessary condition such that the closure of the stratum $S_j$ intersects the stratum $S_i$}, by obvious semicontinuity properties. This condition is {\it  not sufficient}, since for instance the stratum $S_5$ consists of finitely many points, hence its closure coincides with $S_5$, and has empty intersection with the strata  $S_6$ and $S_7$. However, all the other strata $S_j$ contain $g$ in their closure. Note also that if, instead of looking at strata $S_j$ inside $E \cap U_s(2,6)$, we consider the corresponding strata $S'_j$ in $U_s(2,6)$ defined by the same series $HP(H_1(f);t)$, the above claim that $S_5$ is finite it is no longer valid for $S_5'$, and hence it is possible that the closure of $S'_5$ in $U_s(2,6)$ contains $g$. To decide such a question by this approach would involve computation with a family having $\dim M(g)_6=19 $ parameters. On the positive side, the above computation shows that the strata $S'_1$, $S'_2$, $S_3'$, $S_4'$ and $S_6'$ contain $g$ in their closure.
\begin{figure}[ht]
\label{fig:4lines-b}
\setlength{\unitlength}{25pt}
\begin{picture}(6,6)(0,0)
\xygraph{ !{0;<10mm,0mm>:<0mm,10mm>::}
!{(0,0)}*+{S_1}="1"
!{(0,1)}*+{S_2}="2"
!{(-1,2)}*+{S_3}="3"
!{(1,2)}*+{S_4}="4"
!{(0,3)}*+{S_5}="5"
!{(0,4)}*+{S_6}="6"
!{(0,5)}*+{S_7}="7"
"1"-"2" "2"-"3" "2"-"4" "3"-"5" "5"-"4"  "5"-"6" "6"-"7"
}
\end{picture}

\caption{\text{The Hasse diagram of a Hessian  poset in the case $n=2$, $d=6$}}
\label{fig:toy}
\end{figure}
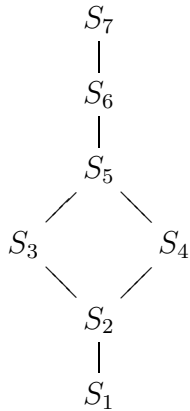

\end{ex}
Here is an example for the smooth curves of degree 6, which can regarded as a continuation of Example \ref{2smooth} and where we consider both series  $HP(H_1(f);t)$ and $HP(H_2(f);t)$. 

\begin{ex}
\label{2smoothC}
Consider the case $n=2$,  $d=6$. Let $g=f''=x^6+y^6+z^6+3x^2y^2z^2$. Then the curve $V(f''):f''=0$ is smooth and we have

\medskip

 \noindent  $HP(H_2(f'');t)=1+3t+6t^2+10t^3+15t^4+18t^5+19t^6+18t^7+9t^8,$
 hence the same series as for the polynomial $f'$ in Example \ref{2smooth} and

\medskip

 \noindent  $HP(H_1(f'');t)=1+3t+6t^2+10t^3+9t^4+6t^5,$
hence a strictly larger series than  for the polynomial $f'$ since the coefficients of $t^5$ satisfy $3<6$.

Consider now the 1-parameter deformation $f_a''=f''+axyz^4$ of $f''$. For small values of $a$ we get

\medskip

 \noindent  $HP(H_2(f''_a);t)=1+3t+6t^2+10t^3+15t^4+18t^5+19t^6+18t^7+9t^8,$
 hence the same series as for the polynomials $f'$ and $f''$. Moreover

\medskip

 \noindent  $HP(H_1(f''_a);t)=1+3t+6t^2+10t^3+9t^4+4t^5,$
hence a  series between that of  the polynomial $f'$ and that of polynomial $f''$, since the corresponding coefficients of $t^5$ satisfy $3<4<6$. However, for the value $a=3$ we get different results, namely

\medskip

 \noindent  $HP(H_2(f''_3);t)=1+3t+6t^2+10t^3+15t^4+18t^5+19t^6+18t^7+9t^8+2t^9,$
  and

\medskip

 \noindent  $HP(H_1(f''_3);t)=1+3t+6t^2+10t^3+9t^4+4t^5,$ hence this polynomial has both series distinct from the minimal values given by the polynomial $f'$ in Example \ref{2smooth}.
Note that the series $HP(H_2(f''_3);t)$ is strictly bigger than the series $HP(H_2(f'');t)$, while the series $HP(H_1(f''_3);t)$ is stricly smaller than $HP(H_1(f'');t)$. {\it The conclusion is that the strata of the Hessian partition to which belong $f''_3$ and respectively $f''$ are not adjacent, i.e. their closures are disjoint}.

\end{ex}

\begin{rk}
\label{rkinflection}
It is not at all clear what the series $HP(H_k(f);t))$ tells us about the geometry of the smooth hypersurface $V(f):f=0$. 

\medskip

\noindent (i) In the curve case, when $n=2$, the equations $f=h(f)=0$ define the {\it inflection points} of the curve $V(f)$, and by Bezout Theorem, their number is at most $d(3d-6)$, with equality for a generic curve, see for instance \cite{CDhess}. This number is lower for some curves, e.g. for the Fermat curve $x^d+y^d+z^d=0$ we have only $3d$ inflection points (counted without multiplicities). However, the series $HP(H_k(f);t))$ do not seem to be related to such differences. Indeed, the two smooth curves  $V(f')$ and $V(f''_3)$ of degree 6 in Example \ref{2smoothC} have each 72 inflection points (by direct computation), but the corresponding series are different.

\medskip

\noindent (ii) The subvariety of $V(f)\subset \PP^n$ given by the vanishing of the ideal $h_{n+2-k}(f)$ corresponds to the closure of the {\it Thom-Boardman singularity set $S_k(\phi_f)$ of the dual mapping} $\phi_f:V(f) \to \PP^n$ associated to $f$. For a generic $f$, the stratum $S_k(\phi_f)$ is smooth, of codimension $k(k+1)/2$ in $V(f)$ (in particular empty if if $n-1 <k(k+1)/2$) and connected if $n-1 >k(k+1)/2$, see \cite{CDhess}.

\medskip

\noindent (iii)  When $d=3$, there is a graded Artinian Gorenstein algebra $R(f)$ with socle  degree 3 associated to $f$ via the so-called {\it Macaulay inverse systems}, see \cite{CRV}, \cite{LP}, \cite{IK}. By results in classical algebraic geometry in \cite{DK}, the variety given by the vanishing of the ideal $h_{k}(f)$ corresponds to the linear forms $\ell \in R(f)_1$ such that the multiplication by $\ell$ gives a linear map $R(f)_1 \to R(f)_2$ of rank $<k$, see Lemma 6.1 in \cite{CRV}. In this way relations to the {\it Lefschetz properties of Artinian graded rings} are established. For such properties, and a  definition of generalized Hessians playing a key role in this theory, see \cite{LP}, in particular Definition 3.75.
\end{rk}


\begin{thebibliography}{00}

\bibitem{CDhess} A.~D.~R.~Choudary, A.~Dimca,  On the dual and hessian mappings of projective 
hypersurfaces, Math. Proc. Cambridge Phil. Soc. 101 (1987), 461-468.


\bibitem{CD} A.~D.~R.~Choudary, A.~Dimca,  Koszul complexes and hypersurface singularities, 
     Proc. Amer. Math. Soc. 121(1994), 1009--1016. 



\bibitem{Co} CoCoA-5 (15 Sept 2014): a system for doing Computations in Commutative Algebra,
available at http://cocoa.dima.unige.it


\bibitem{CRV} A. Conca, M. E. Rossi, G. Valla, Gr\"obner flags and Gorenstein algebras, Compositio Math. 129(2001), 95--121.



\bibitem{Sing} { W. Decker, G.-M. Greuel, G. Pfister \and H. Sch{\"o}nemann.} \newblock {\sc Singular} {4-0-1} --- {A} computer algebra system for polynomial computations, available at {http://www.singular.uni-kl.de} (2014).






\bibitem{D1} A.~Dimca, 
{\em Singularities and topology of hypersurfaces},
Universitext, Springer-Verlag, New York, 1992. 

\bibitem{DBull}  A. Dimca, Syzygies of Jacobian ideals and defects of linear systems,
Bull. Math. Soc. Sci. Math. Roumanie Tome 56(104) No. 2, 2013, 191--203.



\bibitem{DTor}  A. Dimca, Jacobian syzygies, stable reflexive sheaves, and Torelli properties for projective hypersurfaces with isolated singularities, arXiv:1408.2244v3.

\bibitem{DS2} { A. Dimca,  M. Saito},  Generalization of theorems of Griffiths
and Steenbrink to hypersurfaces
with ordinary double points, arXiv:1403.4563v4.

\bibitem{DStNach} A. Dimca, G. Sticlaru, On the syzygies and Alexander polynomials of nodal hypersurfaces, Math.\ Nachr.\ 285 (2012), 2120--2128.


\bibitem{DStCamb} A. Dimca, G. Sticlaru, Chebyshev curves, free resolutions and rational curve arrangements, Math. Proc. Cambridge Phil. Soc. 153  (2012), 385­-397.



\bibitem{DStEdin} A. Dimca, G. Sticlaru, {\em Koszul complexes and pole order filtrations, } Proc. Edinburgh Math. Soc. doi:10.1017/S0013091514000182.

\bibitem{DS14} A. Dimca, E. Sernesi:  {\em Syzygies and logarithmic vector fields along plane curves,} 
Journal de l'\'Ecole polytechnique-Math\'ematiques 1(2014), 247-267.

\bibitem{DK} I. Dolgachev, V. Kanev, Polar covariants of plane cubics and quartics, Adv. Math. 98(1993), 216--301.


\bibitem{EGH} D. Eisenbud, M. Green and J. Harris, Cayley--Bacharach theorems and conjectures,
Bull. Am. Math. Soc. 33 (1996), 295–324.

\bibitem{GLS}
 G.-M. Greuel, C. Lossen and E. Shustin, Equisingular families of projective curves, in Global aspects
of complex geometry, Springer, Berlin, 2006, 171--209.

\bibitem{LP} T. Harima, T. Maeno, H. Morita, Y. Numata, A. Wachi, J. Watanabe, The Lefschetz Properties, Lecture Notes in Math. 2080, Springer, 2013.


\bibitem{Harrris} J. Harris, On the Severi problem, Invent. math. 84 (1986), 445--462.


\bibitem{IK} A. Iarrobino, V. Kanev, Power sums, Gorenstein algebras, and determinantal loci, Lecture Notes in Math. 1721, Springer 1999.





\bibitem{Oka} M. Oka: Symmetric plane curves with nodes and cusps, J. Math. Soc. Japan,
 44 (1992), 375--414.

\bibitem{Se} E. Sernesi: {\it Deformations of Algebraic Schemes},  Springer   Grundlehren b. 334 (2006).

\bibitem{St} G. Sticlaru, Projective Singular Hypersurfaces and Milnor Algebras
A computational approach to classical and new Invariants
LAP LAMBERT Academic Publishing (2014 ).


\end{thebibliography}
\end{document}